\documentclass[11pt]{article}

\usepackage{amsmath}
\usepackage{amssymb}
\usepackage{eucal}
\usepackage{xcolor}

\begin{document}

\baselineskip 16pt

\title{On injectors of Hartley set of a finite group}

\author{Nanying Yang\thanks{Research  is supported by
a NNSF grant of China (Grant \#11301227)and the Natural Science Fundation of Jiangsu Province( grant \#BK20130119).}\\
{\small School of Science, Jiangnan University}\\ {\small Wuxi 214122 , P. R. China}\\
{\small E-mail:
yangny@jiangnan.edu.cn}\\ \\{N.T.Vorob'ev\thanks{Research is supported by the State Research Programme "Convergence" of Belarus (2016-2020).}}\\{\small Department of Mathematics, Masherov Vitebsk State University}\\ {\small Vitebsk 210038, Belarus}\\{\small Email:vorobyovnt@tut.by}\\ \\ {T.B.~Vasilevich\thanks{Research  is supported  by the Belarusian Republican Foundation for Fundamental Research  (F17M-064).}}\\
{\small Department of Mathematics, Masherov Vitebsk State University,} \\ {\small Vitebsk 210038, Belarus} \\
{\small E-mail: tatyana.vasilevich.1992@mail.ru}}

 \date{}
\maketitle

\begin{abstract} Let $G$ be a group and $\mathcal{H}$ be a \emph{Hartley set} of $G$. In this paper, we prove the existence and conjugacy of $\mathcal{H}$-injectors of $G$ and describe the structure of the injectors. As application, some known results are directly followed.

\end{abstract}


\let\thefootnote\thefootnoteorig

\section{Introduction}

Throughout this paper, all groups are finite.  In theory of classes of finite soluble groups, a basic result which generalizes fundamental theorems of Sylow and Hall is the theorem of Fischer, Gasch\"{u}tz and Hartley \cite{fgh} on existence and conjugacy of $\mathfrak{F}$-injectors in soluble groups for every Fitting class $\mathfrak{F}$.

Recall that a class $\mathfrak{F}$ is called a \emph{Fitting class} if $\mathfrak{F}$ is closed under taking normal subgroups and products of normal $\mathfrak{F}$-subgroups. For any class $\mathfrak{F}$ of groups, a subgroup $V$ of a group $G$ is said to be $\mathfrak{F}$-\emph{maximal} if $V\in\mathfrak{F}$ and $U=V$ whenever $V\leq U \leq G$ and $U\in \mathfrak{F}$. From the definition of Fitting class $\mathfrak{F}$, every group $G$ has the largest normal $\mathfrak{F}$-subgroup $G_{\mathfrak{F}}$, so called $\mathfrak{F}$-radical of $G$, which is the product of all normal $\mathfrak{F}$ subgroups.  In particular, if $\mathfrak{F}=\mathfrak{N}$ is the Fitting class of all nilpotent groups, then $G_{\mathfrak{N}}=F(G)$ is the Fitting subgroup of $G$. A subgroup $V$ of a group $G$ is said to be an $\mathfrak{F}$-\emph{injector} of $G$ if $V\cap N$ is an $\mathfrak{F}$-maximal subgroup of $N$ for every subnormal subgroup $N$ of $G$. Note that if $\mathfrak{F}=\mathfrak{N}_{p}$ is the Fitting class of all $p$-groups, then the $\mathfrak{F}$-injectors of a group $G$ are Sylow $p$-subgroups of $G$; if $\mathfrak{F}$ is the Fitting class of all groups with soluble Hall $\pi$-subgroups (i.e $G$ is $E^{s}_{\pi}$-group \cite[p.81]{guo}), where $\pi$ is a set of prime numbers, then the $\mathfrak{F}$-injectors of $G$ are Hall $\pi$-subgroups of $G$.

As a development of the theorem of Fischer, Gasch\"{u}tz and Hartley \cite{fgh}, Shemetkov \cite{shemet1} (resp. Anderson~\cite{ander}) proved that if $G$ is a $\pi$-soluble group (resp. soluble group) and $\mathcal{F}$ is a Fitting set of $G$, then $G$ posseses exactly one conjugacy class of $\mathcal{F}$-injectors, where $\pi$ is the set of all primes divide orders of all subgroups of $G$ in $\mathcal{F}$.

Recall that a nonempty set $\mathcal{F}$ of subgroups of a group $G$ is called a \emph{Fitting set of} $G$ \cite{fgh,shemet1}, if the following three conditions hold:(i) If $T\unlhd S \in \mathcal{F}$ , then $T\in \mathcal{F}$ ; (ii) If $S\in\mathcal{F}$  and $T\in\mathcal{F}$ , $S\unlhd ST$ and $T\unlhd ST$, then $ST\in\mathcal{F}$ ; (iii) If $S\in\mathcal{F}$  and $x\in G$, then $S^{x}\in \mathcal{F}$. So from the definition of Fitting set $\mathcal{F}$, the $\mathcal{F}$-radical $G_{\mathcal{F}}$ of a group $G$ can also be defined as the product of all its normal $\mathcal{F}$-subgroups. For a Fitting set $\mathcal{F}$ of $G$, the $\mathcal{F}$-injector of $G$ is similarly defined as the $\mathfrak{F}$-injector for Fitting class $\mathfrak{F}$ (see [1, Definition~VIII.~(2.5)]).

If $\mathfrak{F}$ is a Fitting class and $G$ is a group, then the set $\{H\leq G: H\in\mathfrak{F}\}$ is a Fitting set, which is denoted by $Tr_{\mathfrak{F}}(G)$ and called the \emph{trace} of $\mathfrak{F}$ in $G$ (see \cite[VIII, 2.2(a)]{Doerk}). Note that for a Fitting class $\mathfrak{F}$, the $\mathfrak{F}$-injectors and $Tr_{\mathfrak{F}}(G)$-injectors of $G$ coincide, but not every Fitting set of $G$ is the trace of a Fitting class (\cite[ VIII. Examples (2.2)(c)]{Doerk}). Hence, if $\mathcal{F}=Tr_{\mathfrak{F}}(G)$, then the theorem of Anderson \cite{ander} and the theorem of Fischer, Gasch\"{u}tz and Hartley \cite{fgh} are Corollaries of the theorem of Shemetkov \cite{shemet1}. Vorob'ev and Semenov \cite{vg} proved that for every set $\pi$ of primes and every Fitting set $\mathcal{F}$ of $\pi$-soluble group $G$, $G$ possesses an $\mathcal{F}$-injector and any two $\mathcal{F}$-injectors are conjugate if $\mathcal{F}$ is $\pi$-saturated, i.e. $\mathcal{F}=\{H\leq G: H/H_{\mathcal{F}}\in \mathfrak{E}_{\pi'}\}$. In connection with these theorems, the following question naturally arise:

\textbf{Question 1.1} \emph{For an arbitrary Fitting set $\mathcal{F}$ of a group $G$ (in case $G$ is a non-$\pi$-soluble group), when $G$ possesses $\mathcal{F}$-injector and any two $\mathcal{F}$-injectors are conjugate?}

There has been substantial research on characterizations of $\mathfrak{F}$-injector for various types of soluble Fitting classes $\mathfrak{F}$ (see \cite{fischer,gv,gl,hartley,hz,mann,seme,shemet1}). It is well known that the product of any two Fitting classes is also Fitting class and multiplication of Fitting classes satisfies associative law (see \cite[Theorem IX.(1.12)(a),(c)]{Doerk}). Hartley \cite{hartley} proved that for the Fitting class of type $\mathfrak{X}\mathfrak{N}$ (where $\mathfrak{X}$ is a nonempty Fitting class and $\mathfrak{N}$ is the Fitting class of all nilpotent groups), a subgroup $V$ of a soluble group $G$ is an $\mathfrak{X}\mathfrak{N}$-injector of $G$ if and only if $V/G_{\mathfrak{X}}$ is a nilpotent subgroup of $G$. As a further improvement, Guo and Vorob'ev \cite{gv} proved that for the Hartley class $\mathfrak{H}$, the set of all $\mathfrak{H}$-injectors of soluble group $G$ coincide with the set of all $\mathfrak{H}$-maximal subgroups of $G$ containing $\mathfrak{H}$-radical of $G$. Let $\mathbb{P}$ be the set of all prime numbers. Following \cite{hartley}, a function $h: \mathbb{P}\rightarrow \{$\emph{nonempty  Fitting classes}$\}$ is a \emph{Hartley function} (or in brevity $H$-\emph{function}). Let $LH(h)=\cap_{p\in\mathbb{P}}h(p)\mathfrak{S}_{p'}\mathfrak{N}_{p}$, where $\mathfrak{N}_{p}$ is the class of all $p$-groups and $\mathfrak{S}_{p'}$ is the class of all soluble $p'$-groups. A Fitting class $\mathfrak{H}$ is called a \emph{Hartley class} if $\mathfrak{H}=LH(h)$ for some $H$-function $h$.

We need to develop and extend the local method of Hartley \cite{hartley}(for soluble Fitting classes) for Fitting sets of groups (not necessary in the universe of soluble groups). For a Fitting set $\mathcal{H}$ of a group $G$ and a nonempty Fitting class $\mathfrak{F}$, we call the set $\{H\leq G: H/H_{\mathcal{H}}\in \mathfrak{F}\}$ of subgroups of $G$ the \emph{product} of $\mathcal{H}$ and $\mathfrak{F}$, and denote it by $\mathcal{H}\circ \mathfrak{F}$ is a Fitting set of $G$ (see Lemma 2.1).

Following \cite{vorobevh}, a function $h: \mathbb{P}\rightarrow \{$\emph{Fitting sets of} $G \}$ is called a Hartley  function of $G$ (or in brevity an $H$-function of $G$).

\textbf{Definition 1.2} Let $h$ be an $H$-function of a group $G$ and $HS(h)=\cap_{p\in\mathbb{P}}h(p)\circ (\mathfrak{E}_{p'}\mathfrak{N}_{p})$, where $\mathfrak{E}_{p'}$ is the class of all $p^{'}$-groups. A Fitting set $\mathcal{H}$ of $G$ called the \emph{Hartley set of} $G$ if $\mathcal{H}=HS(h)$ for some $H$-function $h$.

\textbf{Definition 1.3} Let $\mathcal{H}=HS(h)$ be a Hartley set of a group $G$. Then $h$ is said to be:

(1) \emph{integrated} if $h(p)\subseteq \mathcal{H}$ for all $p$;

(2) \emph{full} if $h(p)\subseteq h(q)\circ \mathfrak{E}_{q'}$ for all different primes $p$ and $q$;

(3) \emph{full integrated} if $h$ is full and integrated as well.

It is easy to see that every Hartley set of a group $G$ can be defined by an integrated $H$-function. Moreover we prove that every Hartley set of $G$ can be defined by a full integrated $H$-function in Lemma 3.4.

In connection with above, the following question naturally arise:

\textbf{Question 1.4 } \emph{Let $G$ be a group (in particular, $G$ is a soluble group), and $\mathcal{H}$ be a Hartley set of $G$, what's the structure of $\mathcal{H}$-injectors of $G$?}

For a Hartley set $\mathcal{H}=HS(h)$ of $G$, $h$ is a full integrated $H$-function of $\mathcal{H}$, we call the subgroup $G_{h}=\prod_{p\in\mathbb{P}}G_{h(p)}$ the $h$-\emph{radical of} $G$. A group $G$ is said to be $\mathfrak{N}$-\emph{constrained} if $C_{G}(F(G))\leq F(G)$. It is well known that if $G$ is soluble, then $G$ is constrained (in general, the converse is not true \cite{iran}).

The following theorem resolved the Questions 1.1 and 1.4.

\textbf{Theorem 1.5 } \emph{Let $\mathcal{H}$ be a Hartley set of a group $G$ defined by a full integrated $H$-function $h$ and $G_{h}$ the $h$-radical of $G$. If $G/G_{h}$ is $\mathfrak{N}$-constrained, then the following statements hold:}

(1) \emph{A subgroup $V$ of $G$ is an $\mathcal{H}$-injector of $G$ if and only if $V/G_{h}$ is a nilpotent injector of $G/G_{h}$;}

(2) \emph{$G$ possesses an $\mathcal{H}$-injector and any two $\mathcal{H}$-injectors are conjugate in $G$;}

(3) \emph{A subgroup $V$ of $G$ is an $\mathcal{H}$-injector of $G$ if and only if $V$ is an $\mathcal{H}$-maximal subgroup of $G$ and $G_{\mathcal{H}}\leq V$.}

Theorem 1.5 give the new theory of $\mathcal{F}$-injectors for Fitting sets of non-soluble groups. From Theorem 1.5, a series of famous results can be directly generalized. For example, Fischer \cite[Corollary IX.4.13]{fischer}, Hartley \cite[section 4.1]{hartley}, Mann \cite[Theorem IX. 4.12]{mann}, Guo and Vorob'ev \cite[Theorem 5.6.8]{guo}.

All unexplained notion and terminology are standard. The reader is referred to \cite{Doerk,guo,ballester}.

\section{Preliminaries}

Note that if all groups in a class $\mathfrak{X}$ are soluble groups (that is $\mathfrak{X}\subseteq\mathfrak{S}$), then $\mathfrak{X}$ is said to be a soluble class.

\textbf{Lemma 2.1} \emph{\cite[Proposition 3.1]{NGV}} \emph{Let $\mathcal{F}$ be a Fitting set of a group $G$ and $\mathfrak{X}$ is a nonempty Fitting class. Then the product $\mathcal{F}\circ \mathfrak{X}$ is a Fitting set of $G$.}

\textbf{Lemma 2.2}  \emph{Let $\mathcal{F}$ and $\mathcal{H}$ be Fitting sets of $G$, and $\mathfrak{X}$, $\mathfrak{Y}$ be nonempty Fitting formations. Then}

(a) \emph{\cite[Proposition 3.4 (3)]{NGV}} $\mathcal{F}\circ (\mathfrak{X}\cap \mathfrak{Y})=\mathcal{F}\circ \mathfrak{X}\cap \mathcal{F}\circ \mathfrak{Y}$.

(b) \emph{\cite[Proposition 3.2 (1)]{NGV}} \emph{If $\mathfrak{M}$ is nonempty Fitting class, then $\mathcal{F}\subseteq\mathcal{F}\circ \mathfrak{M}$.}

(c) \emph{\cite[Proposition 3.4 (2)]{NGV}} $(\mathcal{F}\cap\mathcal{H})\circ\mathfrak{X}=\mathcal{F}\circ\mathfrak{X}\cap \mathcal{H}\circ\mathfrak{X}$.

(d) \emph{\cite[Proposition 3.4 (1)]{NGV}} \emph{If $\mathcal{F}\subseteq \mathcal{H}$, then $\mathcal{F}\circ\mathfrak{X}\subseteq \mathcal{H}\circ\mathfrak{X}$.}

\textbf{Lemma 2.3} \emph{\cite[Proposition 3.3]{NGV}} \emph{Let $\mathcal{F}$ be a Fitting set of a group $G$ and $\mathfrak{X}$, $\mathfrak{Y}$ be Fitting formations. Then $(\mathcal{F}\circ\mathfrak{X})\circ\mathfrak{Y}=\mathcal{F}\circ(\mathfrak{X}\mathfrak{Y})$.}

 \textbf{Lemma 2.4} \emph{\cite[Theorem IV. (1.8)]{Doerk}} \emph{Let $\mathfrak{F}$ and $\mathfrak{H}$ be nonempty formations. If $\mathfrak{F}\subseteq \mathfrak{H}$, then $G^{\mathfrak{H}}\leq G^{\mathfrak{F}}$ for every group $G$.}

\textbf{Lemma 2.5} \emph{\cite[Proposition VIII. (2.4) (d)]{Doerk}} \emph{Let $\mathcal{F}$ be a Fitting set of a group $G$. If $N\unlhd\unlhd G$, then $N_{\mathcal{F}}=N\cap G_{\mathcal{F}}$.}

Let $\mathfrak{F}$ be a nonempty Fitting class. A group $G$ is said to be $\mathfrak{F}$-constrained if $C_{G}(G_{\mathfrak{F}})\leq G_{\mathfrak{F}}$.

\textbf{Lemma 2.6} \emph{\cite[Remark p.624]{Doerk} or \cite{iran}} \emph{The class of all $\mathfrak{N}$-constrained groups is a Fitting class strictly large than $\mathfrak{S}$.}

\textbf{Lemma 2.7} \emph{\cite[Theorem IX.(4.12) (c)-(d)]{Doerk}} \emph{Let $G$ be a group. If $G$ is $\mathfrak{N}$-constrained, then $G$ possesses exactly one conjugacy class of nilpotent injectors.}

The following properties follow at once from definition of an $\mathcal{F}$-injector of a group $G$ and \cite[Remarks IX, (1.3), VIII, (2.6)(2.7)]{Doerk}.

\textbf{Lemma 2.8} \emph{Let $\mathcal{F}$ be a Fitting set of a group $G$ and $\mathfrak{F}$ a class of finite groups. Then}

(a) \emph{If $K\trianglelefteq\trianglelefteq G$ and $V$ is an $\mathcal{F}$-injector of $G$, then $V\cap K$ is an $\mathcal{F}$-injector (or $\mathcal{F}_{K}$-injector) of $K$;}

(b) \emph{If $V$ is an $\mathcal{F}$-injector of $G$, then $G_{\mathcal{F}}\leq V$ and $V$ is an $\mathcal{F}$-maximal subgroup of $G$;}

(c) \emph{If $V$ is an $\mathcal{F}$-maximal of $G$ and $V\cap M$ is an $\mathcal{F}$-injector $M$ for any maximal normal subgroup $M$ of $G$, then $V$ is an $\mathcal{F}$-injector of $G$.}

(d) \emph{If $V \in \mbox{Inj}_{\mathfrak{F}}(G)$ and $\alpha$ : $G \rightarrow G_{\alpha}$ an isomorphism, then $V_{\alpha} \in \mbox{Inj}_{\mathfrak{F}}(G_{\alpha})$; in particular, $\mbox{Inj}_{\mathfrak{F}}(G)$ is a union of $G$-conjugacy classes.}

\textbf{Lemma 2.9} \emph{\cite{mann} or \cite[Theorem IX. (4.12)]{Doerk}} \emph{If $G$ is a $\mathfrak{N}$-constrained group, then a subgroup $V$ of $G$ is a nilpotent injector of $G$ if and only if $F(G)\leq V$ and $V$ is an $\mathfrak{N}$-maximal subgroup of $G$.}

\section{Hartley set and $h$-radical}

In this section we give some results about Hartley sets and $h$-radical of a group $G$, which are also main steps in the proof of Theorem 1.5. Recall that by Lemma 2.1 for a Fitting set $\mathcal{H}$ of $G$ and a nonempty Fitting class $\mathfrak{F}$, the set $\mathcal{H}\circ \mathfrak{F}=\{H\leq G: H/H_{\mathcal{H}}\in\mathfrak{F}\}$ is a Fitting set of $G$. Firstly, we give some following examples of Hartley set.

\textbf{Example 3.1} \emph{(a) Let $\mathcal{N}$ be the trace of the Fitting class $\mathfrak{N}$ in group $G$ and let $h$ be an $H$-function defined as follows: $h(p)=\{1\}$ for all $p\in\mathbb{P}$, where $1$ is an identity subgroup of $G$. Then by Lemma 2.2 (a), we have $HS(h)=\cap_{p\in\mathbb{P}}\{1\}\circ(\mathfrak{E}_{p'}\mathfrak{N}_{p})=\{1\}\circ(\cap_{p\in\mathbb{P}}\mathfrak{E}_{p'}\mathfrak{N}_{p})=\{1\}\circ\mathfrak{N}=\mathcal{N}$. Hence the set of all nilpotent subgroups of $G$ is a Hartley set of $G$.}

\emph{(b) Let $\mathcal{F}$ be a Fitting set of $G$ and $\mathcal{H}=\mathcal{F}\circ\mathfrak{N}$. Let $h$ be an $H$-function such that $h(p)=\mathcal{F}$ for all $p\in\mathbb{P}$. Then by Lemma 2.2 (a) we obtain $HS(h)=\mathcal{F}\circ(\cap_{p\in\mathbb{P}}\mathfrak{E}_{p'}\mathfrak{N}_{p})=\mathcal{F}\circ\mathfrak{N}$ and so $\mathcal{H}$ is a Hartley set of $G$ (the least equality follows from \cite[Lemma 2.7 (a)]{Doerk}).}

\emph{(c) If $k\in\mathbb{N}$, let $\mathcal{N}^{k} (k\geq 1)$ the set of all subgroups of soluble group $G$ of the nilpotent length at most $k$. If $k\geq 1$, take the $H$-function $h$ such that $h(p)=Tr_{\mathfrak{N}^{k-1}}(G)$ for all $p\in \mathbb{P}$. Then by Example (b) above we have $HS(h)=\mathcal{N}^{k}$ is a Hartley set of $G$.}

\emph{(d) Let $\mathcal{H}$ be the trace of Fitting class $\mathfrak{E}_{p'}\mathfrak{N}_{p}$ in group $G$, i.e. $\mathcal{H}$ is the set of all $p$-nilpotent subgroups of $G$. Let $h$ be an $H$-function defined as follows: $h(p)=\{1\}$ and $h(q)=\mathcal{H}$ for all primes $q\neq p$. Then by Lemma 2.2 (a), we have $HS(h)=(\{1\}\circ \mathfrak{E}_{p'}\mathfrak{N}_{p})\cap (\mathcal{H}\circ(\cap_{p\neq q}\mathfrak{E}_{p'}\mathfrak{N}_{p}))=\mathcal{H}\cap \mathcal{H}\circ(\mathfrak{N}_{p}\mathfrak{N}_{p'})$.}

\emph{Now by Lemma 2.2 (b) $HS(h)=\mathcal{H}$ and $\mathcal{H}$ is a Hartley set.}

\textbf{Lemma 3.2} \emph{Every Hartley set can be defined by an integrated $H$-function.}

\textbf{Proof.} Let $\mathcal{H}$ be a Hartley set of a group $G$. Then $\mathcal{H}=HS(h_{1})$ for some $H$-function $h_{1}$. Let $h$ be an $H$-function defined as follows: $h(p)=h_{1}(p)\cap\mathcal{H}$ for all $p\in\mathbb{P}$. By Lemma 2.2 (c), we have

$HS(h)=\cap_{p\in\mathbb{P}}(h_{1}(p)\cap\mathcal{H})\circ\mathfrak{E}_{p'}\mathfrak{N}_{p}=
(\cap_{p\in\mathbb{P}}h_{1}(p)\circ\mathfrak{E}_{p'}\mathfrak{N}_{p})\cap(\cap_{p\in\mathbb{P}}\mathcal{H}\circ(\mathfrak{E}_{p'}\mathfrak{N}_{p}))$.

Hence by Lemma 2.2 (a), (b), $HS(h)=\mathcal{H}\cap \mathcal{H}\circ(\cap_{p\in\mathbb{P}}\mathfrak{E}_{p'}\mathfrak{N}_{p})=\mathcal{H}\cap\mathcal{H}\circ\mathfrak{N}=\mathcal{H}$. The Lemma is proved.

Let $\mathcal{H}$ be a set of subgroups of a group $G$. For $\mathcal{H}$ and nonempty Fitting class $\mathfrak{F}$, we call the set $\{H\leq G: H$ has a normal subgroup $L\in\mathcal{H}$ with $H/L\in\mathfrak{F}$ $\}$ of subgroups of $G$ the product of $\mathcal{H}$ and $\mathfrak{F}$ and denote it by $\mathcal{H}\mathfrak{F}$.

\textbf{Remark 3.3} \emph{Let $\mathcal{H}$ be a Fitting set of a group $G$. Then it is clear that $\mathcal{H}\circ\mathfrak{F}\subseteq\mathcal{H}\mathfrak{F}$. Assume that $\mathfrak{F}$ is the Fitting class and $\mathfrak{F}$ is a homomorph, i.e. $\mathfrak{F}$ is quotient closed. Put $H\leq G$ and $H\in\mathcal{H}\mathfrak{F}$. Then $H/L\in\mathfrak{F}$ for some normal subgroup $L\in\mathcal{H}$ of $H$. Since $L\leq H_{\mathcal{H}}$ and $H/L/H_{\mathcal{H}}/L\cong H/H_{\mathcal{H}}$, $H\in\mathcal{H}\circ\mathfrak{F}$. Thus $\mathcal{H}\mathfrak{F}=\mathcal{H}\circ\mathfrak{F}$.}

Let $G$ be a group and $\mathcal{Y}$ is a set of subgroups of a group $G$. Then Fitset$\mathcal{Y}$ will denote the intersection of all Fitting sets of $G$ that contain $\mathcal{Y}$(see \cite[Definition VIII. 3.1 (b)]{Doerk}).

\textbf{Lemma 3.4} \emph{Every Hartley set can be defined by a full integrated $H$-function.}

\textbf{Proof.} Let $\mathcal{H}$ be a Hartley set of a group $G$. By Lemma 3.2, $\mathcal{H}=HS(h_{1})$ for some integrated $H$-function $h_{1}$. We define a set of subgroups of $G$ by $\overline{h_{1}}(p)=\{ H\leq G: H$ is conjugate with $R^{\mathfrak{S}_{p'}}$ in $G$ for some $R\in h_{1}(p)\}$ for all $p\in\mathbb{P}$. Note that if $H\in \overline{h_{1}}(p)$, then $H\in h_{1}(p)$ and so $\overline{h_{1}}(p)\subseteq  h_{1}(p)$ for all $p\in\mathbb{P}$.

Assume that $X\in \overline{h_{1}}(p)\mathfrak{E}_{p'}$. Then $X$ has a normal subgroup $K\in \overline{h_{1}}(p)$ with $X/K\in \mathfrak{E}_{p'}$. Since $\overline{h_{1}}(p)\subseteq  h_{1}(p)$, $K\leq X_{h_{1}(p)}$. Hence by the isomorphism $X/K/X_{h_{1}(p)}/K\cong X/X_{h_{1}(p)}$, we have $X/X_{h_{1}(p)}\in\mathfrak{E}_{p'}$ and so $X\in h_{1}(p)\circ \mathfrak{E}_{p'}$. Thus $\overline{h_{1}}(p)\mathfrak{E}_{p'}\subseteq h_{1}(p) \mathfrak{E}_{p'}$.

On other hand, let $Y\in h_1(p)\circ \mathfrak{E}_{p'}$. Then $Y/Y_{h_{1}(p)}\in\mathfrak{E}_{p'}$ and $Y^{\mathfrak{E}_{p'}}\leq Y_{h_{1}(p)}$. Hence $Y^{\mathfrak{E}_{p'}}\in h_{1}(p)$. Since $(Y^{\mathfrak{E}_{p'}})^{\mathfrak{E}_{p'}}=Y^{\mathfrak{E}_{p'}}$ and obviously, $Y^{\mathfrak{E}_{p'}}$ is a subgroup conjugate with $Y^{\mathfrak{E}_{p'}}$ in $G$, we have $Y^{\mathfrak{E}_{p'}}\in \overline{h_{1}}(p)$ and so $Y\in \overline{h_{1}}(p)\mathfrak{E}_{p'}$. Thus we obtain the following equation:

$$
\overline{h_{1}}(p)\mathfrak{E}_{p^{'}} = h_{1}(p) \circ \mathfrak{E}_{p^{'}}.                                             \eqno (\ast)
$$

Now, let $h$ be a function such that $h(p) = Fitset(\overline{h_{1}}(p))$ for all $p \in \mathbb{P}$. We prove that $HS(h) = \mathcal{H}$. Since $\overline{h_{1}}(p) \subseteq h_{1}(p)$,  $Fitset(\overline{h_{1}}(p)) \subseteq Fitset(h_{1}(p)) = h_{1}(p)$ and so $h(p) \subseteq h_{1}(p)$.
Hence by Lemma~2.2~(d), we have $h(p) \circ \mathfrak{E}_{p^{'}} \subseteq h_{1}(p) \circ \mathfrak{E}_{p^{'}}$. Note that by Lemma~2.3
$(h(p) \circ \mathfrak{E}_{p^{'}}) \circ \mathfrak{N}_{p} = h(p) \circ (\mathfrak{E}_{p^{'}}\mathfrak{N}_{p})$ and
$(h_{1}(p) \circ \mathfrak{E}_{p^{'}}) \circ \mathfrak{N}_{p} = h_{1}(p) \circ (\mathfrak{E}_{p^{'}}\mathfrak{N}_{p})$.
Therefore by Lemma~2.2~(d) $h(p) \circ (\mathfrak{E}_{p^{'}}\mathfrak{N}_{p}) \subseteq h_{1}(p) \circ (\mathfrak{E}_{p^{'}}\mathfrak{N}_{p})$ for all $p \in \mathbb{P}$. Consequently, $HS(h) \subseteq \mathcal{H}$.

Further more, by ($\ast$), we have $h_{1}(p) \circ \mathfrak{E}_{p^{'}} = Fitset(\overline{h_{1}}(p)\mathfrak{E}_{p^{'}})$. Since $\overline{h_{1}}(p) \subseteq Fitset(\overline{h_{1}}(p))$, by Lemma~2.2~(d) $\overline{h_{1}}(p)\mathfrak{E}_{p^{'}} \subseteq (Fitset(\overline{h_{1}}(p))) \circ \mathfrak{E}_{p^{'}}$. Hence $Fitset(\overline{h_{1}}(p)\mathfrak{E}_{p^{'}}) \subseteq (Fitset(\overline{h_{1}}(p))) \circ \mathfrak{E}_{p^{'}}$.
Thus, for every $p \in \mathbb{P}$, we have the inclusion:

$$
\overline{h_{1}}(p)\mathfrak{E}_{p^{'}} \subseteq (Fitset(\overline{h_{1}}(p))) \circ \mathfrak{E}_{p^{'}} = h(p) \circ \mathfrak{E}_{p^{'}} \eqno  (\ast\ast)
$$

By Lemma 2.3 and Lemma 2.2, we have $h_{1}(p)\circ (\mathfrak{E}_{p'}\mathfrak{N}_{p})\subseteq h(p)\circ(\mathfrak{E}_{p'}\circ\mathfrak{N}_{p})$. Hence $\mathcal{H}\subseteq HS(h)$ and $\mathcal{H}=HS(h)$.

Since $\overline{h_{1}}(p)\subseteq h_{1}(p)$ for all $p\in\mathbb{P}$ and $h_{1}$ is an integrated $H$-function of $\mathcal{H}$, $h(p)\subseteq\mathcal{H}$ for all $p\in\mathbb{P}$ and so $h$ is an integrated $H$-function of $\mathcal{H}$.

Now, we show that $h(p)\subseteq h(q)\circ \mathfrak{E}_{q'}$ for all $p\neq q$.

Let $H$ be an arbitrary subgroup in $h_{1}(p)$ and $p\neq q$. Since $h_{1}$ is an integrated $H$-function, $H\in\mathcal{H}$ and so $H\in h_{1}(q)\circ(\mathfrak{E}_{q'}\mathfrak{N}_{q})$. By Lemma 2.3 $h_{1}(q)\circ(\mathfrak{E}_{q'}\mathfrak{N}_{q})=(h_{1}(q)\circ\mathfrak{E}_{q'})\circ\mathfrak{N}_{q}$. Hence $H^{\mathfrak{N}_{q}}\in h_{1}(q)\circ\mathfrak{E}_{q'}$. Since $p\neq q$, $H^{\mathfrak{E}_{p'}}\leq H^{\mathfrak{N}_{q}}$ by Lemma 2.4. Consequently, $H^{\mathfrak{E}_{p'}}\in h_{1}(q)\mathfrak{E}_{q'}$. Then by ($\ast \ast$), $H^{\mathfrak{E}_{p'}}\in h(q)\circ \mathfrak{E}_{q'}$ for every $H\in h_{1}(p)$ and all primes $p\neq q$.

Let $R\in \overline{h_{1}}(p)$. Then by definition of the set $\overline{h_{1}}(p)$, we have that $R$ is a conjugate subgroup of $G$ with $S^{\mathfrak{E}_{p'}}$ for some subgroup $S\in h_{1}(p)$. Therefore $R\in h(q)\mathfrak{E}_{q'}$ and so $h_{1}(p)\subseteq h(q)\circ\mathfrak{E}_{q'}$. Thus $h(p)=Fitset(\overline{h_{1}}(p))\subseteq Fitset(h(q)\circ \mathfrak{E}_{q'})=h(q)\circ \mathfrak{E}_{q'}$ for all primes $q\neq p$. This complete the proof.

Recall that $G_{h}=\prod_{p\in\mathbb{P}}G_{h(p)}$, where $h$ is a full integrated $H$-function of a Hartley set of a group $G$, i.e. $h$ is an $h$-radical of $G$.

\textbf{Lemma 3.5} \emph{Let $\mathcal{H}$ be a Hartley set of a group $G$ and $h$ is a full integrated $H$-function of $\mathcal{H}$. If $H$ is a subgroup of $G$ such that $G_{h}\leq H$ and $H/G_{h}$ is a nilpotent subgroup of $G/G_{h}$, then $H\in\mathcal{H}$.}

\textbf{Proof.} Let $q$ be an arbitrary prime number. Since $G_{h(q)}\unlhd H$, $G_{h(q)}\leq H_{h(q)}$ for all $q\in\mathbb{P}$. Let $p\in \mathbb{P}$ and $p\neq q$. Note that $G_{h(q)}G_{h(p)}/G_{h(q)}\cong G_{h(p)}/G_{h(q)}\cap G_{h(p)}=G_{h(p)}/(G_{h(p)})_{h(q)}$. Since $h$ is a full integrated $H$-function of $\mathcal{H}$, $h(p)\subseteq h(q)\mathfrak{E}_{q'}$. Hence $G_{h(p)}\in h(q)\mathfrak{E}_{q'}$. Therefore $G_{h(q)}G_{h(p)}/G_{h(q)}\in\mathfrak{E}_{q'}$ for all primes $p$ and $q$. Consequently, $G_{h}/G_{h(q)}\in\mathfrak{E}_{q'}$ and by using the isomorphisms $H_{h(q)}G_{h}/H_{h(q)}\cong G_{h}/H_{h(q)}\cap G_{h}\cong (G_{h}/G_{h(q)})/((H_{h(q)}\cap G_{h})/G_{h(q)})$, we obtain that $H_{h(q)}G_{h}/H_{h(q)}$ is a $q'$-group. Since $H/G_{h}$ is a nilpotent subgroup of $G/G_{h}$, $H/H_{h(q)}G_{h}$ is also a nilpotent subgroup of $G/G_{h}$ and so $H/H_{h(q)}G_{h}\in \cap_{p\in\mathbb{P}}\mathfrak{E}_{q'}\mathfrak{N}_{q}$. Therefore, by the isomorphism $H/H_{h(q)}G_{h}\cong (H/H_{h(q)})/(H_{h(q)}G_{h}/H_{h(q)})$, we have that $H/H_{h(q)}\in\mathfrak{E}_{q'}\mathfrak{N}_{q}$ for all $q\in\mathbb{P}$. Hence $H\in\cap_{q\in\mathbb{P}}h(q)\mathfrak{E}_{q'}\mathfrak{N}_{q}=\mathcal{H}$. The Lemma is proved.

\textbf{Lemma 3.6} \emph{Let $\mathcal{H}$ be a Hartley set of a group $G$. If $h$ is a full integrated $H$-function of $\mathcal{H}$, then $G_{\mathcal{H}}/G_{h}=F(G/G_{h})$.}

\textbf{Proof.} Let $F(G/G_{h})=R/G_{h}$. Since $h$ is a integrated $H$-function of $\mathcal{H}$, $(G_{\mathcal{H}})_{h(p)}=G_{h(p)}$. Hence $G_{\mathcal{H}}/G_{h(p)}\in \mathfrak{E}_{p'}\mathfrak{N}_{p}$ for all $p\in\mathbb{P}$ and so $G_{\mathcal{H}}/G_{h}$ is a nilpotent subgroup of $G/G_{h}$. Hence $G_{\mathcal{H}}/G_{h}\leq F(G/G_{h})$ and we have that $G_{\mathcal{H}}\leq R$.

On the other hand, since $R/G_{h}$ is a nilpotent subgroup of $G/G_{h}$, by Lemma 3.5, $R\in \mathcal{H}$. Hence $R\leq G_{\mathcal{H}}$. Thus $F(G/G_{h})=G_{\mathcal{H}}/G_{h}$. The Lemma is proved.

\textbf{Lemma 3.7} \emph{Let $h$ be a full integrated $H$-function of a Hartley set $\mathcal{H}$ of a group $G$. If $G/G_{h}$ is $\mathfrak{N}$-constrained, $G_{\mathcal{H}}\leq H\leq G$ and $H\in\mathcal{H}$, then $H/G_{h}$ is nilpotent.}

\textbf{Proof.} Since $G_{\mathcal{H}}\unlhd H$ and $h$ is an integrated $H$-function of $\mathcal{H}$, by Lemma 2.5, $G_{h(p)}=(G_{\mathcal{H}})_{h(p)}=H_{h(p)}\cap G_{\mathcal{H}}$. Hence $[H_{h(p)},G_{\mathcal{H}}]\leq H_{h(p)}\cap G_{\mathcal{H}}=G_{h(p)}$ and so $H_{h(p)}\leq C_{G}(G_{\mathcal{H}}/G_{h(p)})\leq C_{G}(G_{\mathcal{H}}/G_{h})$. Since $G/G_{h}$ is constrained and by Lemma 3.6 $G_{\mathcal{H}}/G_{h}=F(G/G_{h})$, $C_{G}(G_{\mathcal{H}}/G_{h})\leq G_{\mathcal{H}}$. Thus $H_{h(p)}\leq G_{\mathcal{H}}$. Therefore, $G_{h(p)}=(G_{\mathcal{H}})_{h(p)}=H_{h(p)}\cap G_{\mathcal{H}}=H_{h(p)}$ for all $p\in \mathbb{P}$. Hence, $H/G_{h}=H/H_{h}\in \mathfrak{E}_{p'}\mathfrak{N}_{p}$ for all $p\in\mathbb{P}$ and so $H/G_{h}$ is a nilpotent group. This complete the proof.

\textbf{Corollary 3.8} \emph{Let $h$ be a full integrated $H$-function of a Hartley set of a group $G$. Let $G/G_{h}$ be an $\mathfrak{N}$-constrained group and $G_{\mathcal{H}}\leq H\leq G$. Then $H\in\mathcal{H}$ if and only if $H/G_{h}$ is nilpotent.}

\section{Proof and Some Applications of Theorem 1.5}

\textbf{Proof of Theorem 1.5} (1) We first prove that if $V$ is a subgroup of a group $G$ such that $V/G_{h}$ is a nilpotent injector of $G/G_{h}$, then $V$ is an $\mathcal{H}$-injector of $G$. We prove this statement by induction on the order of $G$.

Let $M$ be an arbitrary maximal normal subgroup of $G$ and $M_{h}$ is an $h$-radical of $M$.

Since $h$ is a full integrated $H$-function of Hartley set $\mathcal{H}$, $h(p)\subseteq h(q)\circ\mathfrak{E}_{q'}$ and $h(p)\subseteq \mathcal{H}$, for all different $p, q \in \mathbb{P}$. Then, by the isomorphism $G_{h(q)}G_{h(p)}/G_{h(q)}\cong G_{h(p)}/G_{h(p)}\cap G_{h(q)}=G_{h(p)}/(G_{h(p)})_{h(q)}$, we see that $G_{h(q)}G_{h(p)}/G_{h(q)}$ is a $q'$-group of $G/G_{h(q)}$ for any prime $q$. Hence $G_{h}/G_{h(q)}$ is also a $q'$-group. Since $(G_{h}\cap M)G_{h(q)}/G_{h(q)}\leq G_{h}/G_{h(q)}$, by the isomorphism $(G_{h}\cap M)G_{h(q)}/G_{h(q)}\cong (G_{h}\cap M)/(G_{h}\cap M)\cap G_{h(q)}=(G_{h}\cap M)/G_{h(q)}\cap M$, by Lemma 2.5, we obtain that $G_{h}\cap M/M_{h(q)}$ is a $q'$-group, for all $q\in\mathbb{P}$. Now, note that $(G_{h}\cap M)/M_{h(q)}/M_{h}/M_{h(q)}\cong G_{h}\cap M/M_{h}$. Hence $G_{h}\cap M/M_{h}\in \cap_{q\in\mathbb{P}}\mathfrak{E}_{q'}=(1)$ and so $G_{h}\cap M=M_{h}$.

If $G_{h}\leq M$. Then $G_{h}=M_{h}$.

Since $V/G_{h}$ is a nilpotent injector of $G/G_{h}$, $V \cap M/G_{h}$ is a nilpotent injector of $M/G_{h}$ and consequently $V \cap M/M_{h}$ is a nilpotent injector of $M/M_{h}$. Hence, by induction, $V \cap M$ is an $\mathcal{H}$-injector of $M$.

Now, in order to complete the proof of the statement, we only need to prove that $V$ is an $\mathcal{H}$-maximal subgroup of $G$. Since $V/G_{h}$ is nilpotent and $G_{\mathcal{H}}\leq V $, by Lemma 3.5 $V\in\mathcal{H}$. Assume that $V< V_{1}$, where $V_{1}$ is an $\mathcal{H}$-maximal subgroup of $G$. Since $V\cap M$ is an $\mathcal{H}$-maximal subgroup of $M$, $V\cap M=V_{1}\cap M$. Hence $V_{1}$ is an $\mathcal{H}$-maximal subgroup of $G$ and $V \cap M$ is an $\mathcal{H}$-injector of $M$ for any maximal normal subgroup $M$ of $G$. Consequently, by Lemma 2.8 (b), (c), $V_{1}$ is an $\mathcal{H}$-injector of $G$ and $G_{\mathcal{H}}\leq V_{1}$. Then, by Corollary~3.8, we have $V_{1}/G_{h}$ is a nilpotent subgroup of $G/G_{h}$, contrary to the fact that $V/G_{h}$ is $\mathfrak{N}$-maximal in $G/G_{h}$. Hence $V=V_{1}$ and so by Lemma 2.8 (c), $V$ is an $\mathcal{H}$-injector of $G$.

If $G_{h}\nleq M$. In this case, by the maximality of $M$, we have that $G = G_{h}M$. Since, by Lemma~2.8(d), $G/G_{h}\cong M/G_{h}\cap M = M/M_{h}$ and so $V\cap M/M_{h}$ is a nilpotent injector of $M/M_{h}$. Then, by induction, $V\cap M$ is an nilpotent injector of $M$. By Lemma~3.5, we know that $V\in\mathcal{H}$. If $V < F_{1}$, where $F_{1}$ is an $\mathcal{H}$-maximal subgroup of $G$, then $V\cap M = F_{1}\cap M$. Since $G_{h}\leq V$, $VM = G$. Consequently, $F_{1} =  F_{1}\cap VM = V(F_{1}\cap M) = V(V\cap M) = V$ and so $V$ is an $\mathcal{H}$-maximal subgroup of $G$. Therefore $V$ is an $\mathcal{H}$-injector of $G$.

Conversely, let $V$ be an $\mathcal{H}$-injector of $G$. We prove that $V/G_{h}$ is a nilpotent injector of $G/G_{h}$. By Lemma 2.8 (b),(c) $G_{\mathcal{H}} \leq V$ and $V$ is $\mathcal{H}$-maximal in $G$. Hence by Lemma 3.5, $V/G_{h}$ is nilpotent. Since $V$ is $\mathcal{H}$-maximal in $G$, we obtain that $V/G_{h}$ is $\mathfrak{N}$-maximal subgroup of $G/G_{h}$ containing the nilpotent radical of $G/G_{h}$. Consequently, by Lemma 2.9, $V/G_{h}$ is a nilpotent injector of $G/G_{h}$. Thus, statement (1) hold.

(2) The existence of $\mathcal{H}$-injectors has proved in (1). Let $F/G_{h}$ and $V/G_{h}$ are nilpotent injectors of $G/G_{h}$. Then by Lemma 2.7 $F/G_{h}$ and $V/G_{h}$ are conjugate in $G/G_{h}$. Hence $F$ and $V$ are conjugate in $G$.

(3) Let $V$ be an $\mathcal{H}$-injector of $G$. Then by Lemma 2.8 (b), (c), $V\geq G_{\mathcal{H}}$ and $V$ is $\mathcal{H}$-maximal in $G$.

Conversely, let $V$ be an $\mathcal{H}$-maximal subgroup of $G$ and $V\geq G_{h}$. We prove that $V$ is an $\mathcal{H}$-injector of $G$. Clearly, $G_{h}\leq V$. Then by Lemma 3.5, $V/G_{h}$ is nilpotent. Since $V$ is $\mathcal{H}$-maximal in $G$, $V/G_{h}$ is $\mathfrak{N}$-maximal in $G/G_{h}$. Now by Lemma 3.6, $V/G_{h}\geq F(G/G_{h})$. Hence by Lemma 2.9, $V/G_{h}$ is a nilpotent $\mathcal{H}$-injector of $G/G_{h}$. This complete the proof of the theorem.

Let $\mathcal{X}$ be a Fitting set of a group $G$ and $\mathcal{X}\circ \mathfrak{S}=\{H\leq G: H/H_{\mathcal{X}}\in \mathfrak{S}\}$, where $\mathfrak{S}$ is a class of all soluble groups. Note that the set $\mathcal{X}\circ\mathfrak{S}$ is a Fitting set by Lemma 2.1.

We give some applications of our main results

\textbf{Corollary 4.1} \emph{Let $G\in\mathcal{X}\circ\mathfrak{S}$ and $\mathcal{H}=\mathcal{X}\circ \mathcal{N}$ is a Fitting set of $G$. Then:}

\emph{(1) A subgroup $V$ of $G$ is an $\mathcal{H}$-injector of $G$ if and only if $V/G_{\mathcal{X}}$ is nilpotent injector of $G/G_{\mathcal{X}}$.}

\emph{(2) The set of all $\mathcal{H}$-injectors of $G$ is exactly the subgroups $V$ of $G$ such that $V\geq G_{\mathcal{H}}$ and $V$ is $\mathcal{H}$-maximal in $G$.}

\textbf{Proof.} By Example 3.1 (b), $\mathcal{H}$ is a Hartley set of $G$ , which can be defined by full integrated $H$-function $h$ such that $h(p)=\mathcal{X}$ for all $p\in\mathbb{P}$. Since $G/G_{\mathcal{X}}$ is soluble, $G/G_{\mathcal{X}}$ is $\mathfrak{N}$-constrained.

Let $\mathfrak{X}$ be a nonempty Fitting class and $\mathfrak{H}=\mathfrak{X}\mathfrak{N}$ is a Fitting product of $\mathfrak{X}$ and $\mathfrak{N}$, then we have the following result immediately from our Theorem.

\textbf{Corollary 4.2} \emph{(1) A subgroup $V$ of a group $G\in \mathfrak{F}\mathfrak{S}$ is an $\mathfrak{H}$-injector of $G$ if and only if $V/G_{\mathfrak{X}}$ is a nilpotent injector of $G/G_{\mathfrak{X}}$.}

\emph{(2) A subgroup $V$ of a group $G\in \mathfrak{F}\mathfrak{S}$ is an $\mathfrak{H}$-injector if and only if $V\geq G_{\mathfrak{H}}$ and $V$ is $\mathfrak{H}$-maximal in $G$.}

\textbf{Proof.} Since the set of all $\mathfrak{X}\mathfrak{N}$-injectors of $G$ and the set of all $Tr_{\mathfrak{X}\circ\mathfrak{N}}(G)$-injectors of $G$ are coincide, so Corollary 4.2 holds from Corollary 4.1.

\textbf{Corollary 4.3} \emph{(Hartley \cite{hartley})} \emph{Let $\mathfrak{X}$ be a nonempty soluble Fitting class and $\mathfrak{H}=\mathfrak{X}\mathfrak{N}$. A subgroup $V$ of soluble group $G$ is an $\mathfrak{H}$-injector if and only if $V/G_{\mathfrak{X}}$ is a nilpotent injector of $G/G_{\mathfrak{X}}$.}

\textbf{Corollary 4.4} \emph{(Fischer \cite{fischer})} \emph{A subgroup $V$ of soluble group $G$ is a nilpotent injector of $G$ if and only if $F(G)\leq V$ and $V$ is $\mathfrak{N}$-maximal in $G$.}

\textbf{Corollary 4.5} \emph{(Guo and Vorob'ev \cite{gv})} {Let $\mathfrak{H}$ be a soluble Hartley class and $G$ is a soluble group. Then a subgroup $V$ of $G$ is an $\mathfrak{H}$-injector of $G$ if and only if $V/G_{h}$ is a nilpotent injector of $G/G_{h}$.}

\textbf{Corollary 4.6} \emph{Let $\mathcal{N}^{k} (k\geq 1)$ be a Fitting set of all subgroup of a soluble group $G$ with a nilpotent length at most $k$. Then the set of all $\mathcal{N}^{k}$-injectors of $G$ is exactly the set of all subgroups $V$ of $G$ such that $V/G_{\mathcal{N}^{k-1}}$ is a nilpotent of $G/G_{\mathcal{N}^{k-1}}$. In particular, a subgroup $V$ is an $\mathcal{N}^{2}$-injector of $G$ if and only if $V/F(G)$ is a nilpotent injector of $G/F(G)$.}

\textbf{Proof.} By Example 3.1 (c), $\mathcal{N}^{k}$ is a Hartley set of $G$ and the function $h$ such that $h(p)=Tr_{\mathfrak{N}^{k-1}}(G)$ for all $p\in \mathbb{P}$ is a full integrated $H$-function of $\mathcal{N}^{k}$. Since a group $G$ is soluble, $G/G_{h}$ is $\mathfrak{N}$-constrained and so the Corollary hold from our Theorem.

\textbf{Corollary 4.7} \emph{(Guo and Vorob'ev \cite{gv})} \emph{Let $\mathfrak{N}^{k}$ $(k\geq 1)$ be the class of all groups with nilpotent length at most $k$ and $G$ a soluble group. Then the set of all $\mathfrak{N}^{k}$-injectors of $G$ is exactly the set of all subgroups $V$ of $G$ such that $V/G_{\mathfrak{N}^{k-1}}$ is a nilpotent injector of $G/G_{\mathfrak{N}^{k-1}}$. In particular, a subgroup $V$ of a soluble group $G$ is metanilpotent injector if and only if $V/F(G)$ is a nilpotent injector of $G/F(G)$.}

\

\end{document}